\documentclass[12pt]{amsart}
\usepackage{a4wide,times,amssymb}

\newif\ifpdf
\ifx\pdfoutput\undefined
  \pdffalse
\else
  \pdfoutput=1
  \pdftrue
\fi

\ifpdf
  \usepackage[usenames,pdftex]{color}
  \usepackage[pdftex]{hyperref}
\else
  \usepackage[usenames,dvips]{color}
  \usepackage[hypertex]{hyperref}
\fi

\title{Marcinkiewicz-Zygmund inequalities}
\author{Joaquim Ortega-Cerd\`a}
\address{Dept.\ Matem\`atica Aplicada i An\`alisi, Universitat  de Barcelona,
Gran Via 585, 08007 Bar\-ce\-lo\-na, Spain}
\email{jortega@ub.edu}
\thanks{The authors are supported by the DGICYT grant: BFM2002-04072-C02-01
and by the CIRIT grant: 2001SGR00172}
\author{Jordi Saludes}
\address{Dept. Matem\`atica Aplicada 2, Universitat Polit\`ecnica de Catalunya,
Edifici TR5, Colom 11, 08222 Terrassa, Spain}
\email{jordi.saludes@upc.edu}
\date{\today}
\newcommand{\D}{\mathbb D}
\newcommand{\C}{\mathbb C}
\newcommand{\R}{\mathbb R}
\newcommand{\Q}{\mathbb Q}
\newcommand{\N}{\mathbb N}
\newcommand{\T}{\mathbb T}
\newcommand{\ZZ}{\mathcal Z}
\newcommand{\WW}{\mathcal W}
\renewcommand{\P}{\mathcal P}
\newtheorem{theorem}{Theorem}
\newtheorem{prop}{Proposition}
\newtheorem{lemma}{Lemma}
\newtheorem{corollary}{Corollary}
\theoremstyle{definition}
\newtheorem*{example}{Example}
\newtheorem*{definition}{Definition}

\begin{document}

\begin{abstract}
We study a generalization of the classical Marcinkiewicz-Zygmund inequalities. 
We relate this problem to the sampling sequences in the Paley-Wiener
space and by using this analogy we give sharp necessary and sufficient
computable conditions for a family of points to satisfy the
Marcinkiewicz-Zygmund inequalities.
\end{abstract}
\maketitle
\section{Introduction} We recap the classical Marcinkiewicz-Zygmund
inequalities (see \cite{MarZyg37} or \cite[Theorem~7.5, chapter~X]{Zygmund77}).
Let $\omega_{n, j}$, $j = 0, \ldots, n$ be
the $( n + 1 )$-roots
of the unity. We denote by $\P_n$ the holomorphic polynomials of degree
smaller or equal than $n$.
Then for any $q\in\P_n$ we have
\begin{equation}
  \label{M-Z} \frac{C_p^{- 1}}{n} \sum_{j = 0}^n |q ( w_{n, j} ) |^p \le
  \int_0^{2 \pi} |q ( e^{i \theta} ) |^p \, d \theta \le
  \frac{C_p}{n} \sum_{j = 0}^n |q ( w_{n, j} ) |^p,
\end{equation}
for any $1 < p < + \infty$. The essential feature is that $C_p$ is independent
of the polynomial $q$ and of \textit{the degree} of the polynomial. We aim at
generalizing these inequalities to more general families of points.

We will consider a triangular family of points $z_{n, j} \in \T$ of the form
\[
\ZZ = \{ z_{n, j} \}_{\begin{subarray}{c}
\\n = 0, \ldots, \infty\\j = 0, \ldots, m_n \end{subarray}}.
\]
We will denote by $\ZZ(n)$ the $n$-th generation of points in the family, i.e.
$\ZZ(n)=\{z_{n,0},\ldots,z_{n,m_n}\}$.
\begin{definition}
  We say that $\ZZ$ is a M-Z (Marcinkiewicz-Zygmund) family for $L^p$ ($1 \le p
  < + \infty$) if the following inequality holds for all holomorphic polynomials
$q$ of
  degree smaller or equal than $n$
  \[ \frac{C_p^{- 1}}{m_n} \sum_{j = 0}^{m_n} |q ( z_{n, j} ) |^p \le
     \int_0^{2 \pi} |q ( e^{i \theta} ) |^p \, d \theta \le
     \frac{C_p}{m_n} \sum_{j = 0}^{m_n} |q ( z_{n, j} ) |^p . \]
  Of course $m_n \ge n$ for all $n$. When $p = \infty$ the inequality is
  replaced by
  \[ \sup_{|z| = 1} |q ( z ) | \le C \sup_{j = 0, \ldots, m_n} |q ( z_{n, j} )
     |. \]
\end{definition}

This sort of inequalities are similar to the sampling sequences in the
Paley-Wiener setting. We will show that this similitude is more than
superficial and show how the same kind of results are expected.

\textit{A minimal M-Z family of points} is a M-Z family such that $m_n = n$.
These have been studied and described in detail by Chui and Zhong in
\cite{ChuZho99} when $1 < p < \infty$. If $p = 1$ or $p = \infty$ there are
no minimal M-Z families (see Theorem~\ref{nohiha}) but there are plenty of M-Z
families. When $1 < p < \infty$ a naive guess suggests that any M-Z family of
points minus some points maybe a  minimal M-Z family. The
following
example shows that this is not the case and one cannot reduce the study of M-Z
families to the minimal
ones.

\begin{example}
  Take the triangular family $\ZZ = \{ z_{n, j} = e^{2\pi ij / ( n + 2 )}
  \}_{\begin{subarray}{c}j = 0, \ldots, ( n + 1 ) \\
  n = 0, \ldots, \infty\end{subarray}}$. Clearly $\ZZ$ is a M-Z family 
  for $L^p$ ($1<p<\infty$) but there is not
  any triangular subfamily $\WW$ that  is a minimal M-Z family.
\end{example}

\begin{proof}
  We will argue by contradiction. Assume that $\WW\subset \ZZ$ is a minimal
subfamily in $L^p$. 
  In each generation $n$ of $\ZZ$ there is an excess of one point. Since the
  problem is
  invariant under rotations we may assume that the minimal family $\WW$ is
  just $\ZZ$ minus the point $1$ in all generations. Consider the polynomials
  $p_{n} ( z ) = 1 + z +\cdots +z^n$. The norm of $p_{n}$ can be
  easily estimated with the classical M-Z inequality,
  \[ \| p_{n} \|^p_{L^p} = \int_{|z| = 1} |p_{n} ( z ) |^p\, |dz| \simeq
     \frac{|p_{n} ( 1 ) |^p}{n+1}   = (n+1)^{p-1}, \text{ since }
     p_n(z)=\frac{z^{n+1}-1}{z-1},\ \text{for }z\ne 1. \]
  On the other hand if $\WW$ is a M-Z family then
  \[ \| p_{n} \|^p_{L^p} \simeq \frac{1}{n+1} \sum_{j = 1}^{n+1} |p_{n} (
     e^{2\pi ij / ( n + 2 )} ) |^p, \]
  but
  \[ |p_{n} ( e^{2\pi ij / ( n + 2 )} ) | = \left| \frac{1 - e^{-2\pi ij  / (
n+2
     )}}{1 - e^{2\pi ij / ( n + 2 )}} \right| =1, \]
  which yields a contradiction.
\end{proof}

There are several possible motivations for this work. One possible motivation
is the approximation of periodic continuous functions by
trigonometric polynomials. Consider for
instance any triangular family of points $\WW$ such that $\WW(n)$ has
cardinality $2n+1$. There are periodic continuous functions $f$ such that the 
unique trigonometric polynomial of degree $n$ that interpolates $f$ in $\WW(n)$
does not converge (in uniform norm) to $f$ (see \cite{Cheney98}). To obtain a
convergent sequence of trigonometric polinomials $p_n$ to $f$ it is possible to
use the following proposition:
\begin{prop}
Let $f\in C(\T)$ 
and let  $\ZZ$ be a M-Z family for $L^\infty$. If $p_n$ is the
trigonometric polinomial of degree $n$ that minimizes 
$\max_{z\in \ZZ(2n)} |p_n(z)-f(z)|$ then $p_n\to f$ in $L^\infty(\T)$.
\end{prop}
\begin{proof}
First observe that if $\ZZ$ is a M-Z family for the holomorphic polynomials
with the norm $L^p(\T)$ then $\WW$ defined as
$\WW(n)=\ZZ(2n)$ is an $L^p$ M-Z family for the harmonic polynomials with
the norm $L^p(\T)$. The reason is that for any harmonic polynomial $\pi$ of
degree $n$ ($\pi(z)=a_0+\sum_{1\le i\le n} a_i z^i+b_i \bar z^i$), the
polynomial $p=z^n\pi$, when restricted to $\T$, coincides with a
holomorphic polynomial of degree $2n$. Moreover $|p(e^{ix})|=|\pi(e^{ix})|$ for
all $x\in \R$, thus the $L^p$ norm of $\pi$ and $p$ are the same and the
discretized norms are the same too. Therefore the description of M-Z families
for harmonic polynomials can be reduced to the study of M-Z families of
holomorphic polynomials. We will, as usual, identify any periodic function on
$\R$ with a function in $\T$ and the trigonometric polynomials with the harmonic
polynomials. 

Take the function $f\in C(\T)$. There exists a  sequence of harmonic polynomials
$q_n$ of degree $n$ that converge to $f$ in the uniform norm by Weierstrass
Theorem. Let $p_n$ be the harmonic polynomials of degree $n$ that minimize
$\max_{z\in \WW(n)} |p_n(z)-f(z)|$. Clearly $\|p_n -f\|_\infty\le
\|p_n-q_n\|_\infty
+\|q_n-f\|_\infty$ and $\|q_n-f\|_\infty\to 0$. Moreover since $\WW$ is a M-Z
family
$\|p_n-q_n\|_\infty\lesssim \max_{z\in \WW(n)} |p_n(z)-q_n(z)|$. Since $p_n$
minimizes the distance to $f$ in $\WW(n)$ we have then that 
$\|p_n-q_n\|_\infty\lesssim \max_{z\in \WW(n)} |q_n(z)-f(z)|\to 0$. Thus
$p_n\to f$ in
the uniform norm. 
\end{proof}
A full characterization of the M-Z families for $p=\infty$ is given by
Theorem~\ref{six}.

There is  also
some motivation in the study of M-Z families that comes from a problem in
Computerizedd Tomography. In the setting of the Radon transform in dimension
two,
one typically knows the integrals of a function supported in the
unit disk through a finite number of lines, and one wants to reconstruct the
function from the value of these integrals. The lines (in the usual
parallel-beam geometry)
are grouped in families of parallel lines along a finite number of directions.
The number of directions (which can be identified with points in the unit
circle) depends on the resolution that we want to achieve. Typically it is
required 
that the set of directions is a uniqueness set for the polynomials of a
certain degree. But it has been noted (see \cite[p.~70]{Natterer86} or
\cite[p.~668]{Logan75})
that uniqueness is not good enough for the numerical stability of the
reconstruction. We need certain stability conditions. This is exactly what
the M-Z inequalities provide. Thus, in principle, the M-Z family of
points provide good sets of directions to sample the Radon transform. A more
detailed analysis of the application of M-Z families to the Computerized
Tomography deserves a work of its own. 

In light of this connection it seems also
interesting to study analogous M-Z inequalities  in higher dimensions (i.e.,
replace the circle by the sphere in $\R^3$ and the holomorphic
polynomials by harmonic polynomials of a certain degree). Some preliminary work
has already been done, \cite{MhaNarWar01} but we don't pursue this line further.

We will rather provide metric conditions for  $\ZZ$ to be a M-Z family. Our
first main result is Theorem~\ref{tfive} which gives a sharp metric condition
for a
family $\ZZ$ to be M-Z. This condition is in terms of a density. When
$p=\infty$ the density condition is actually a characterization. This is our
other main result (Theorem~\ref{six}). As mentioned before all this results
are parallel to similar results for entire functions in the Paley-Wiener space.
A good reference for these is \cite{Seip04}. In the next section we prove this
metric characterization after some preliminary technical lemmas. Finally in the
last section we briefly comment on a full characterization of M-Z families when
$p=2$. This characterization is in terms of the invertibility of certain
Toeplitz operators and it is somehow involved. We have not been able to
obtain good computable conditions nor interesting examples from it.

\section{Metric conditions}

We state now some preliminary results that we will need for our computation. The
following inequality was found by Bernstein and Zygmund, see
for instance \cite[Theorem 3.16, p. 11 vol II]{Zygmund77}.

\begin{theorem}
  [Bernstein type inequalities]\label{bernstein}For any $p$, $0 < p \le
  \infty$, and any polynomial $q_n$ of degree $n$:
  \[ \| q_n' \|_{L^p ( \mathbb{T} )} \le n \| q_n \|_{L^p ( \mathbb{T} )} \]
\end{theorem}

There is a good reason that the classical M-Z inequality does not hold in the
endpoints cases $p=1,\infty$. It is not true in this case, but in the irregular
setting that we consider this is still the case.
\begin{theorem}\label{nohiha}
  There are no Marcinkiewicz-Zygmund families for $p = 1$ and $m_n = n$.
\end{theorem}

\begin{proof}
  Suppose that there exists such family $\ZZ$. We start by proving that in
any generation $\ZZ(n)$ two different points $z_{n,j}$ and $z_{n,k}$, $j\ne
k$ are uniformly separated. More precisely there is a constant $C>0$ such that
$n|z_{n,j}-z_{n,k}|\ge C$ for all $j\ne k$ and all $n$. To prove this, take
the unique polynomial $p\in\P_n$ with values $p(z_{n,j})=1$, $p(z_{n,k})=0$
for $k=0,\ldots, n$, $k\ne j$. Since we assume that $\ZZ$ is M-Z,
then $n\|p\|_1\simeq 1$. The Bernstein inequalities entails that
$\|p'\|_1\le n\|p\|_1\simeq 1$ and the classical M-Z inequality
 implies that $\|p'\|_\infty\le n\|p'\|_1\simeq n$. Thus 
\[
1=|p(z_{n,j})-p(z_{n,k})|\le \|p'\|_\infty |z_{n,j}-z_{n,k}|\lesssim n
|z_{n,j}-z_{n,k}|.
\]
Now we will  build a bounded
  projection from $L^1 ( \T )$ to the Hardy space $H^1 ( \T )$ and this is well
  known to be impossible, thus we will reach a contradiction. To build such
  projection, take any function $f \in L^1 ( \T )$
  and for any $n$, consider the values of the Poisson extension $v_j = P [ f ]
  ( z_{n, j} ( 1 - 1 / n ) )$, $j = 0, \ldots, n$. We denote by $p(z,e^{it})$
the Poisson kernel in the disk. An easy but tedious computation shows that
$\sup_{t}\sum_{j=0}^n p((1-1/n)z_{n,j},e^{it})\le C n$, because the points in
$\ZZ(n)$ satisfy $n|z_{n,j}-z_{n,k}|\ge C'$. Thus,
$|v_0 | + \cdots +   |v_n | \lesssim n \| f \|_1$. If $\ZZ$ is a M-Z family then
the holomorphic
  polynomial $p_n$ that takes the values $v_j$ at the points $z_{n, j} ( 1 - 1
  / n )$ has norm $\| p_n \|_1$ bounded by $\| f \|_1$. Therefore the operator
  $Q_n$ that associates at each function $f$ the corresponding
interpolating polynomial
  $p_n$ is a bounded projection from $L^1$ to the subspace of holomorphic
  polynomials of degree $n$. Now take a a partial subsequence of $Q_n$
  converging to $Q$. This is a bounded projection into $H^1$ (the polynomials
  are dense and are fixed by $Q$).
\end{proof}

Observe that the same proof shows that there are no minimal M-Z families for
$p = \infty$.

It will convenient to evaluate the norm of a polynomial not on the boundary of
the unit disk, but on the boundary of some slightly smaller or bigger disk.
This can be done without harm as Lemma~\ref{lema1} shows.

We need the following lemma:
\begin{lemma}[Hardy]\label{hardy} Let $p>0$,
  let $f$ be holomorphic on $D ( 0, R )$ and define
  \[ I ( r ) = \int_0^{2 \pi} |f ( re^{i \theta} ) |^p \, d \theta . \]
  Then $I ( r )$ is increasing and log-convex with respect to $\log r$. 
\end{lemma}
This is a classical result of Hardy~\cite{Hardy15} that actually follows from
elementary theory of subarhomic functions (see e.g. \cite[Theorem
2.16]{HayKen76}).
%
%

\begin{lemma}\label{perturbatiu}
  \label{lema1}Let $p \in [ 1, \infty ]$ and let $q$ be any polynomial of
  degree $n$. For any $r \in [ \frac{n}{n+1}, \frac{n+1}n ]$ there is a constant
  $C_p$ (independent of $n$ and $q$) such that
  \begin{equation}
    C_p \| q \|_{L^p} \le \| q_r \|_{L^p} \le C_p^{- 1} \| q \|_{L^p},
    \label{eq:l1}
  \end{equation}
  where $q_r$ is the dilation $q_r ( z ) = q ( rz )$.
\end{lemma}

\begin{proof}
  Let $q \in \P_n$. 
  In the case $0<r < 1$ and $p<\infty$  
  we obtain $\| q_r \|_p \leq \| q \|_p$ because
  $\| q_r \|_p^p$ is increasing by Lemma~\ref{hardy}. When $p = \infty$, the
   conclusion  follows from the maximum principle since $|q(z)|^p$ is
subharmonic.
  
  For $r > 1$, define
  \[ I_p ( r ) = \int_0^{2 \pi} |q ( re^{i \theta} ) |^{p} d \theta
      \text{ and }  I_{\infty} ( r ) = \max_{\theta \in
     [ 0, 2 \pi ]} |q ( re^{i \theta} ) |. \]
  By using Hadamard's three-circle principle  for $p = \infty$ and
  Lemma~\ref{hardy} for $p <\infty $ , we can assume that $I_p$ is
log-convex as a function of $\log r$ for
  $p \in [ 1, \infty ]$. Therefore
  \begin{equation}
    \log I_p ( r ) \leq ( 1 - t ) \log I_p ( 1 ) + t \log I_p ( R )
    \label{eq:h3c}
  \end{equation}
  where $t = \log r / \log R$. Notice that $I_p ( r ) = \| q_r \|_p$.
  
  Now: $I_p ( R ) = O ( R^n )$ as $R
\rightarrow
  \infty$, thus (\ref{eq:h3c}) becomes:
  \[
    \log I_p ( r ) \leq \log I_p ( 1 ) + \varepsilon ( R ) + n \log r
  \]
  with $\varepsilon ( R ) \rightarrow 0$; so we have
  \begin{equation}
    \| q_r \|_p \le \| q \|_p r^n\le e\|q\|_p, \text{ for }1<r<\frac{n+1}n
\label{eq:qrn} .
  \end{equation} 
  
  For the left hand side of \eqref{eq:l1} in the case $1 < r$, we can use
  again Lemma~\ref{hardy} for $|q ( z )|^{p}$, in the case $p<\infty$ and the
  subharmonicity of the absolute value of $q$ in the case $p = \infty$.
  
  In the case $1 - 1 / n < r < 1$, we consider $\tilde{q} ( z ) = q ( rz )
  \in \P_n$, now $\tilde{q}_{r^{- 1}} ( z ) = q ( z )$ with $1 < r^{- 1} < n
  / ( n - 1 )$. Therefore, using \eqref{eq:qrn} we have
  \[ \| q \|_p = \| \tilde{q}_{r^{- 1}} \|_p \leq C \| \tilde{q} \|_p = C \|
     q_{r^{}} \|_p . \]
  
\end{proof}

If we denote by $C_n$ the annulus $\{ z \in \mathbb{C} : 1 - 1 / n < |z| < 1 +
1 / n \}$ and $dm ( z )$ the Lebesgue measure, Lemma~\ref{perturbatiu}
immediately entails the following corollary

\begin{corollary}\label{unic}
  For any polynomial $q$ of degree $n$
  \[ \| q \|_{L^p ( \mathbb{T} )}^p \simeq n \| q \|_{L^p ( C_n, dm ( z ) )}^p .
  \]
\end{corollary}

Now we are able to prove a Plancherel-Polya type Theorem describing the
triangular families that satisfy the first of the M-Z
inequalities (the easier one)

\begin{theorem}
  \label{PlPo}Let $p \in [ 1, \infty )$. If $\mathcal{Z}$ is a triangular
  family such that
  \begin{equation}
    \label{rd} \# ( \mathcal{Z} ( n ) \cap I_n ) \frac{n}{m_n} \le C
  \end{equation}
  for all $n \in \mathbb{N}$ and all intervals $I_n$ of the unit circle of
  length $1 / n$, then for any polynomial $q$ of degree $n$
  \begin{equation}
    \label{carleson} \frac{1}{m_n} \sum_{k = 0}^{m_n} |q ( z_{n, k} ) |^p \le
    C_p \int_0^{2 \pi} |q ( e^{it} ) |^p \, dt,
  \end{equation}
  where the constant $C_p$ is independent of the degree. Conversely if
  \eqref{carleson} holds, then there is a constant $C$ such that \eqref{rd}
holds
  for all intervals $I_n$ of length $1 / n$.
\end{theorem}

\begin{proof}
  Take any point $z_{n, k} \in \mathcal{Z} n$. By the subharmonicity of $|q|^p$
  we have
  \[ |q ( z_{n, k} ) |^p \le \frac{n^2}{\pi} \int_{D ( z_{n, k}, 1 / n )} |q (
     w ) |^p \, dm ( w ) . \]
  Now if we add all the points we get
  \[ \frac{1}{m_n} \sum_{k=0}^{m_n} |q ( z_{n, k} ) |^p \le \frac{n^2}{\pi m_n}
\sum_{k=0}^{m_n}
     \int_{D ( z_{n, k}, 1 / n )} |q ( w ) |^p \, dm ( w ) . \]
  and now we replace the sum in the right hand side by the integral over the
  union of disks. Each point in the annulus $C_n$ is at most in $C
m_n/n$ disks due to the hypothesis~\ref{rd}. Finally the sum is bounded by:
  \[ \frac{1}{m_n} \sum |q ( z_{n, k} ) |^p \le \frac{n^2}{\pi m_n}
     \frac{Cm_n}{n} \int_{||z| - 1| < 1 / n} |q ( w ) |^p \, dm (
     w ) . \]
  Finally we can apply Corollary~\ref{unic}
  \[ n \int_{||z| - 1| < 1 / n} |q ( w ) |^p \, dm ( w ) \simeq
     \| q \|^p_{L^p(\T)} . \]
  From now on, we will use the following notation for the discrete norm:
  \[ \| q \mid \mathcal{Z} ( n ) \|_p^p = \frac{1}{m_n} \sum_{k = 0}^{m_n} |q
     ( z_{n, k} ) |^p, \quad \text{for } q \in \P_n \]
  where $\mathcal{Z} ( n ) = \{ z_{K, 0}, \ldots, z_{K, m_N} \}$ is the $n$-th
  slice of the triangular family $\mathcal{Z}$ and $m_n = \# \mathcal{Z} ( n
  )$.
  
  For the second part, consider the polynomial
  \[ q_m ( z ) = \frac{z^m - 1}{m ( z - 1 )} = \frac{1 + z + z^2 + \cdots +
     z^{m - 1}}{m} \]
  This polynomial satisfies $\| q_m \|_{\infty} = 1$ and moreover $q_m ( 1 ) =
  1$. Let $\mathcal{W} = \{ w_{m, j} \}$ be the triangular family of the
  $m$-roots
  of the unity ($w_{m, j} = e^{i 2 \pi j / m}$, for $j = 0, \ldots, m - 1$).
  We have $q_m ( w_{m, j} ) = 0$ for $j \neq 0$ and $q_m ( w_{m, 0} ) =
  1$. If we fix $p \geq 1$, it is clear that $\| q_m \|_p \simeq \| q_m \mid
  \mathcal{W} ( m ) \|_p = m^{- 1}$ because the roots of unity are the
  prototypical M-Z family. So
  \begin{equation}
    \| q_m \|_{p^{}}^p \leq C_p ( \mathcal{W} ) \| q_m \mid \mathcal{W} ( m )
    \|_p^p = C_p ( \mathcal{W} ) / m \label{mz-w}
  \end{equation}
  Now assume that \eqref{rd} is false, but \eqref{carleson} is true for a
  given constant $C_p$. Then, by taking $C = 2^{p + 2} C_p ( \mathcal{W} )
  C_p$ in the reverse of \eqref{rd}, there is $N > 0$ and an arc $I$ of length
  $1 / N$ such that
  \[ \# ( \mathcal{Z} ( N ) \cap I ) > Cm_N / N. \]
  Now, divide $I$ in halves; it is clear that there is a half $J$ such that
  \begin{equation}
    \# ( \mathcal{Z} ( N ) \cap J ) \geq Cm_N / ( 2 N ) . \label{half}
  \end{equation}
  Since $\| q ( e^{i \theta} \cdot ) \|_p = \| q \|_p$ and $\| q ( e^{i
  \theta} \cdot ) \mid \mathcal{Z} ( n ) \|_p = \| q \mid e^{i \theta}
  \mathcal{Z} ( n ) \|_p$, we can assume that $J$ is centered at the point 1,
  without changing the MZ property of $p_N$. On the other hand, by the
  Bernstein inequality,
  \[ \sup_{|z| = 1} |q'_N ( z ) | \le N \sup_{|z| = 1} |q_N ( z ) | = N ; \]
  Combined with $|q_N ( z ) - q_N ( 1 ) | \leqslant \sup_{| \xi | = 1} |
  \nabla q_N ( \xi ) ||z - 1|$ and the fact that $|J| = ( 2 N )^{- 1}$, we
  obtain a lower bound for $q_N$ on $z \in J$:
  \[ |q_N ( z ) | \geq 1 - N|z - 1| \geq 1 / 2. \]
  Then, using \eqref{half} in the definition of the discrete norm
  \begin{equation}
    \| q_N \mid \mathcal{Z} ( N ) \|_p^p \geq \frac{1}{m_N} \inf_{z \in J}
    |q_N ( z ) |_p^p \# ( \mathcal{Z} ( N ) \cap J ) \geq \frac{1}{m_N}
    \frac{1}{2^p} C \frac{m_N}{2 N} \geq 2 C_p C_p ( \mathcal{W} ) / N
    \label{nomz1}
  \end{equation}
  But now, by \eqref{mz-w} and the assumption that $\mathcal{Z}$ satisfies
  \eqref{carleson},
  \begin{equation}
    \| q_N \mid \mathcal{Z} ( N ) \|_p^p \leq C_p \| q_N \|_p^p \leq C_p C_p (
    \mathcal{W} ) / N \label{nomz2}
  \end{equation}
  Inequalities \eqref{nomz1} and \eqref{nomz2} are incompatible, thus
  \eqref{carleson} cannot hold for $\mathcal{Z}$.
\end{proof}

\begin{definition}
  Given a triangular family $\mathcal{Z}$ we say that it is separated whenever
  there is an $\varepsilon > 0$ such that $|z_{n, j} - z_{n, k} | \ge
  \varepsilon / n$ for all $1 \le j, k \le m_n$, $j \neq k$ and all $n \in
  \mathbb{N}$.
\end{definition}

\begin{theorem}
  If $\mathcal{Z}$ is a M-Z family then there is a separated subfamily
  $\mathcal{Z}'$ such that $\mathcal{Z}'$ is also a M-Z family.
\end{theorem}

In view of this theorem we will limit ourselves to the study of separated
triangular families. Observe that any separated triangular family satisfies
$m_n \le Cn$.

\begin{proof}
  The idea of the proof is the following. We take $\varepsilon > 0$ very
  small (to be determined) and we split the circle $|z| = 1$ into intervals
  $I_n$ of size $\varepsilon / n$. From the points belonging to $\mathcal{Z} (
  n )$ we are going to select some to be in $\mathcal{Z}' ( n )$. In each
  interval $I_n$ we only keep at most one point. If the remaining points are
  still not $\varepsilon / ( 3 n )$-distance one from the other we discard
  some more points in such a way that all points in $\mathcal{Z}' ( n )$ are
  at least $\varepsilon / ( 3 n )$-distance one from the other and any point
  in $\mathcal{Z} ( n )$ is at most at distance $3 \varepsilon / n$ from some
  of the points in $\mathcal{Z}' ( n )$. We need now to prove that
  $\mathcal{Z}'$ is a M-Z family for a small enough $\varepsilon > 0$.
  
  To begin with need the following stability result.
  
  \begin{lemma}
    If $\mathcal{Z}$ is a M-Z triangular family then there is an $\epsilon >
    0$ (depending only on the constants of the M-Z inequalities for
    $\mathcal{Z}$) such that for any perturbation $\mathcal{Z}^{\ast}$ of the
    original family with the property $|z_{n, j} - z_{n, j}^{\ast} | \le
    \varepsilon / n$ is still a M-Z triangular family. 
  \end{lemma}
  
  \begin{proof}
    Observe that
    \[ \left| \Bigl(\frac{1}{m_n} \sum_{j = 1}^{m_n} |q_n ( z_{n, j} ) |^p\Bigr
       )^{1/p}-
       \Bigl(\frac{1}{m_n} \sum_{j = 1}^{m_n} |q_n ( z_{n, j}^{\ast} ) |^p
        \Bigr)^{1/p}\right|
       \le \Bigl(\frac{1}{m_n} \sum_{j = 1}^{m_n} |q_n ( z_{n, j}^{\ast} ) - q_n
(
       z_{n, j} ) |^p\Bigr)^{1/p} . \]
    There are points $\tilde{z}_{n, j}$ in between $z_{n, j}^{\ast}$ and
    $z_{n, j}$ such that
    \[ |q_n ( z_{n, j}^{\ast} ) - q_n ( z_{n, j} ) |^p \le C_p |q_n' (
       \tilde{z}_{n, j} ) |^p |z_{n, j} - z_{n, j}^{\ast} |^p \]
    Now
    \[ \frac{1}{m_m} \sum_{j = 1}^{m_n} |q_n ( z_{n, j}^{\ast} ) - q_n ( z_{n,
       j} ) |^p \le \frac{C_p \varepsilon^p}{n^p m_n} \sum_{j = 1}^{m_n} |q_n' (
       \tilde{z}_{n, j} ) |^p \]
    The points in the triangular family $\tilde{z}_{n, j}$ satisfy \eqref{rd}
    because $z_{n, j}$ does and they are very close one to the other.
    Therefore we can apply Theorem~\ref{PlPo} and we get
    \[ \frac 1{m_n}\sum_{j = 1}^{m_n} |q_n ( z_{n, j}^{\ast} ) - q_n ( z_{n, j}
) |^p \le
       \frac{C_p\varepsilon^p}{n^p} \int_{\mathbb{T}} |q_n' |^p \, dt. \]
    Finally we can use Bernstein inequalities (Theorem~\ref{bernstein}) and
    the fact that $z_{n, j}$ is a M-Z family and we get
    \[ \left|\Bigl( \frac{1}{m_n} \sum_{j = 1}^{m_n} |q_n ( z_{n, j} ) |^p\Bigr
)^{1/p}-
       \Bigl(\frac{1}{m_n} \sum_{j = 1}^{m_n} |q_n ( z_{n, j}^{\ast} )
|^p\Bigr)^{1/p} \right|
       \le \frac{1}4\Bigl(\frac 1{ m_n} \sum_{j = 1}^{m_n} |q_n ( z_{n, j} )
|^p\Bigr)^{1/p} \]
    if we pick $\varepsilon$ small enough. Therefore
    \[ \frac{1}{m_n} \sum_{j = 1}^{m_n} |q_n ( z_{n, j} ) |^p \simeq
       \frac{1}{m_n} \sum_{j = 1}^{m_n} |q_n ( z_{n, j}^{\ast} ) |^p \]
    as we wanted to prove.
  \end{proof}
  
  We finish now the proof of the theorem. Since the family $\mathcal{Z}'$ is
  $\varepsilon / 3$ separated we have automatically the inequality
  \eqref{carleson}. We only have to prove the other inequality.
  
  For any point $z_{n, j} \in \mathcal{Z}_n$ the closest point $z_{n,
  j}^{\ast}$ in $\mathcal{Z}_n$ is at most at distance $3 \varepsilon / n$, 
  so we can apply the Lemma. We can't conclude directly that
  $\mathcal{Z}'$ is a M-Z family because in the discrete norm we may be
  repeating the same $z_{n, j}^{\ast}$ associated to many different $z_{n,
  j}$. The inequality \eqref{rd} does the trick: there is a bound of at most
  $C \frac{m_n}{n}$ different $z_{n, j}$ points in $\mathcal{Z} ( n )$
  associated to the same point $z_{n, j}^{\ast} \in \mathcal{Z}' ( n )$.
  \[ \| q_n \|^p \simeq \frac{1}{m_n} \sum_{i = 1}^{m_n} |q_n ( z_{n,
     j}^{\ast} ) |^p \le \frac{1}{m_n} \sum_{i = 1}^{m_n'} \frac{Cm_n}{n} |q_n
     ( z_{n, j}' ) |^p . \]
  Since $\mathcal{Z}'$ is separated then $m_n' \simeq n$ and thus
  \[ \| q_n \|^p \simeq \frac{1}{m_n'} \sum_{i = 1}^{m_n'} |q_n ( z_{n, j}' )
     |^p \]
\end{proof}

In the statements of Theorem~\ref{tfive} and Theorem~\ref{six} we denote by 
$(x,y)$ to the arc in $\T$ delimited by the endpoints $e^{i x}$ and $e^{i y}$.
\begin{theorem}\label{tfive}
  Given a separated family $\mathcal{Z}$, if
  \[ D^- ( \mathcal{Z} ) = \liminf_{R \to \infty} \left( \liminf_{n \to
     \infty} \frac{\min_{x \in [0,2\pi]} \# \mathcal{Z} ( n ) \cap ( x, x +
     R / n )}{R} \right) > \frac 1{2 \pi}, \]
  then $\mathcal{Z}$ is a M-Z family (for any $p \in [ 1, \infty ]$).
  Conversely, if $\mathcal{Z}$ is a M-Z family for some $p \in [ 1, \infty ]$,
  then
  \begin{equation}
    \label{necessaria} D^- ( \mathcal{Z} ) = \liminf_{R \to \infty} \left(
    \liminf_{n \to \infty} \frac{\min_{x \in [0,2\pi]} \# \mathcal{Z} ( n )
    \cap ( x, x + R / n )}{R} \right) \ge \frac 1{2 \pi} .
  \end{equation}
\end{theorem}

\begin{proof}
  We start with the sufficiency part for $p = \infty$. We will relate this
  problem to the similar problem in the Bernstein class which consists on
  entire functions of exponential type $\pi$ bounded on the real line. The
  sampling sequences for such functions were studied and described by Beurling
  in \cite[p. 340]{Beurling89b}.
  
  Let $\mathcal{Z}$ be a separated triangular family. To each generation of
  points
  \[ \mathcal{Z} ( n ) = \{ e^{i \theta_{n, 1}}, e^{i \theta_{n, 2}}, \ldots,
     e^{i \theta_{n, m_n}} \}, \quad \theta_{n, i} \in [ - \pi, \pi ], \]
  we associate a real sequence $\Lambda ( n )$ consisting of the points
  \begin{equation}
    \label{aplanats} \Lambda ( n ) = \{ n \theta_{n, 1} / ( 2 \pi ) + nk, n
    \theta_{n, 2} / ( 2 \pi ) + nk, \ldots, n \theta_{n, m_n} / ( 2 \pi ) + nk
    \}_{k \in \mathbb{Z}} .
  \end{equation}
  Since $\mathcal{Z}$ is separated then $\Lambda ( n )$ is $\delta$-separated
  uniformly on $n$. Moreover the hypothesis on $\mathcal{Z}$ imply that there
  is an $R > 0$ and $\varepsilon > 0$ such that
  \[ \frac{\# \Lambda ( n ) \cap ( x, x + R )}{R} > 1 + \varepsilon, \quad
     \forall x \in \mathbb{R} . \]
  This means that $\Lambda ( n )$ is a sampling sequence for the Bernstein
  class (see \cite[p. 346 Theorem 5]{Beurling89b}. That is
there is a
  constant $C$ which depends only on $\varepsilon, R$ and the separation
  constant $\delta$ such that $\sup_{\mathbb{R}} |f ( x ) | \le C
  \sup_{\lambda \in \Lambda ( n )} |f ( \lambda ) |$ for all functions $f$ in
  the Bernstein class. The constant is independent of $n$.
  
  Given any polynomial $q \in \mathcal{P}_n$ we have $\sup_{w \in \mathbb{T}}
  |q ( w ) | = \sup_{x \in \mathbb{R}} |q ( e^{2 \pi ix / n} ) |$. If we
  define $f \in \mathcal{H} ( \mathbb{C} )$ as $f ( w ) = q ( e^{2 \pi iw / n}
  ) e^{- \pi iw}$ then $f$ belongs to the Bernstein class since $q$ is of
  degree $n$. Therefore we may apply Beurling's Theorem and we obtain
  \[ \| q \|_{L^{\infty} ( \mathbb{T} )} = \| f \|_{L^{\infty} ( \mathbb{R} )}
     \le C \sup_{\lambda \in \Lambda ( n )} |f ( \lambda ) | = \sup_{z_i \in
     \mathcal{Z} ( n )} |q ( z_i ) |. \]
  Thus we have proved the theorem for $p = \infty$. Now we are going to prove
  it for $p = 1$ and the others will follow by interpolation. We will use a
  similar scheme as in \cite[p. 36]{Seip93} Indeed, the property that
  $\mathcal{Z}$ is a M-Z family for $p$ means that the operators $R_n : (
  \mathcal{P}_n, \| \cdot \|_p ) \to ( \mathbb{C}^{m_n}, \| \cdot \|_p )$
  defined as $R_n ( q ) = ( q ( z_{n, 1} ), \ldots, q ( z_{n, m_n} ) )$ are
  injective and of closed range. Therefore the inverse $R_n$ is defined in the
  range of $R_n$ and it has bounded norm $\| R_n^{- 1} \|_p$. The key point is
  that the norm of the inverse must be bounded by $Cn^{- 1 / p}$. We have
  proved that whenever $D^- ( \mathcal{Z} ) > 1$ then $\| R_n^{- 1}
  \|_{\infty} < C$. We will now prove that $\| R_n^{- 1} \|_1 < C / n$ is also
  uniformly bounded, and by interpolation $\| R_n^{- 1} \|_p < Cn^{- 1 / p}$
  for any $p \in [ 1, \infty ]$.
  
  We will use that $\mathcal{Z}$ is a M-Z family for $p = \infty$. Let us
  denote by $( A_n, \| \cdot \|_{\infty} ) \subset \mathbb{C}^{m_n}$ the image
  of $R_n$. Any bounded linear functional $\phi$ on $( \mathcal{P}_n, \|
  \|_{\infty} )$ induces a bounded linear functional $\phi$ on $A_n$ as
  $\tilde{\phi} ( x ) = \phi ( R^{- 1} ( x ) )$, with $\| \tilde{\phi} \| \le
  K \| \phi \|$. For each $w \in \mathbb{T}$ let $\phi_w$ denote the point
  evaluation functional,i.e, $\phi_w ( q ) = q ( w )$ for any $q \in
  \mathcal{P}_n$. The norm of $\phi$ is trivially $1$. Since the dual space of
  $(\mathbb{C}^{m_n}, \| \cdot \|_{\infty})$ is $(\mathbb{C}^{m_n}, \| \cdot
  \|_1)$, there is a $m_n$-tuple of numbers $g_j ( w )$ such that $\sum_{j =
  1}^{m_n} |g_j ( w ) | \le M$ and moreover
  \begin{equation}
    \label{reconstruccio} q ( w ) = \sum_{j = 0}^{m_n} q ( z_{n, j} ) g_j ( w
    ) .
  \end{equation}
  Moreover since there is an $\epsilon > 0$, such that $D^- ( \mathcal{Z} ) >
  1 + \varepsilon$ then the $MZ$-inequality holds not only for polynomials $q$
  of degree $n$ but also on polynomials of degree $[( 1 + \varepsilon / 2 ) n]$.
  Thus we have established \eqref{reconstruccio} for all polynomials $q$ of
  degree $( 1 + \varepsilon / 2 ) n$. Consider now a collection of auxiliary
  polynomials $a_n ( z )$ of at most degree $[ \varepsilon n / 2 ]$ such that
  $a_n ( 1 ) = 1$ and $\| a_n \|_1 \simeq 1 / n$. This polynomial can be
  constructed for instance taking $a_n ( z ) = b_n^2 ( z )$ and $b_n ( z )$ 
  a polynomial of degree $[ \varepsilon n / 8 ]$ which is $1$ in $1$ and $0$
  in the other $[ \varepsilon n / 8 ]$-roots of unity. Clearly since the roots
  of unity are a M-Z family $\| b_n \|_2^2 \simeq \frac{1}{n}$. Moreover $\|
  a_n \|_1 = \| b_n \|_2^2$. Finally, take any polynomial $r$ of degree $n$
  and any point $w \in \mathbb{T}$. The polynomial $q ( z ) = r ( z ) a_n (
  \bar{w} z )$ is a polynomial of degree at most $( 1 + \varepsilon / 2 ) n$
  with the property that $q ( w ) = r ( w )$. We may apply \eqref{reconstruccio}
  and we get
  \[ r ( w ) = \sum_{j = 0}^{m_n} r ( z_{n, j} ) a_n ( \bar{w} z_{n, j} ) g_j
     ( w ) . \]
  If we now estimate $\| r \|_1$ we get
  \[ \| r \|_1 \le \sum_{j = 0}^{m_n} |r ( z_{n, j} ) | \sup_j
     \int_{\mathbb{T}} |a_n ( \bar{w} z_{n, j} ) g_j ( w ) | \,
     d|w|. \]
  But $|g_j ( w ) | \le M$ (even the sum is bounded by $M$) and $\int_{|w| =
  1} |a_n ( \bar{w} z_{n, j} ) | \, d|w| = \| a_n \|_1 \simeq 1 /
  n$, therefore
  \[ \| r \|_1 \lesssim \frac{1}{n} \sum_{j = 0}^{m_n} |r ( z_{n, j} ) |, \]
  for all polynomials $r$ of degree $n$ which is what we wanted to prove.
  
  To prove the necessity we want to deal only with $p=2$. The next lemma shows 
  how we can reduce ourselves to this situation.
  
  \begin{lemma}
    \label{p-estabilitat}If $\ZZ$ is a separated $L^p$ M-Z family then for any
    arbitrary small $\delta > 0$ the family $\ZZ'$ obtained scaling the
    indexes (i.e. $\ZZ' ( n ) = \ZZ ( [ n ( 1 + \delta ) ] )$ is an $L^2$ M-Z
    family.
  \end{lemma}
  
  \begin{proof}
    We will prove that under the hypothesis $\ZZ'$ is a M-Z family for $L^1$
    and for $L^\infty$, thus by interpolation it will be a M-Z family for all
    $L^r$, $1\le r\le \infty$, in particular for $r=2$ as in the statement. 
    We start by proving that $\ZZ'$ is an $L^\infty$ M-Z family.
    Just as before if $\ZZ$ is a $L^p$ M-Z family then there are functions
    $g_j^n : \T \to \C$ such that $\sum_{j=0}^{m_n} |g_{n,j} ( z ) |^q \le C$
    (where $q$ satisfies $1/p+1/q=1$) and for all polynomials of degree $n$:
    \[ p ( z ) = \sum^{m_n}_{j = 0} p ( z_{n, j} ) g_{n,j} ( z ) . \]
    
    If we take polynomials $c_n$ of degree $[ \delta n ]$ such
    that $\| c_n \|_p \simeq n^{-1/p}$ and $c_n(1)=1$,  we get that for any $z
     \in \T$,
     \begin{equation*}\label{pessos} 
     p ( z ) = \sum_{j = 0}^{m_n} p ( z'_{n, j} ) c_n ( z'_{n, j} \bar{z} )
       g_{n,j} ( z ), 
      \end{equation*}
      for all polynomials of degree $n$ and the rescaled sequence $\ZZ'$.
     If we use H\"older inequality we obtain
     \[ |p(z)|\le (\sup_j |p(z'_{n,j}|) \|c_n(z'_{n, j} \bar{z})\|_{\ell^p}
     \|g_j^n(z)\|_{\ell^q}.\]
      Finally, by Theorem~\ref{PlPo} $\|c_n(z'_{n, j}
      \bar{z})\|_{\ell^p}\lesssim      n^{1/p}\|c_n\|_{p}$ and thus
      \[
      \sup_{\T} |p(z)|\lesssim  (\sup_j |p(z'_{n,j})|).
      \]
      That proves that $\ZZ'$ is an $L^\infty$ M-Z family. To prove that it is 
      an $L^1$ family we take polynomials $b_n$ of degree $[ \delta n ]$ such
    that $\| a_n \|_1 \simeq n^{-1}$ and $a_n(1)=1$,  we get that for any $z
     \in \T$,
     \[
     p ( z ) = \sum_{j = 0}^{m_n} p ( z'_{n, j} ) a_n ( z'_{n, j} \bar{z} )
       g_{n,j} ( z ), 
      \]
      if we integrate this 
      \[
      \|p\|_1\le \sum_{j = 0}^{m_n}\int_\T |p ( z'_{n, j} )| |a_n (z'_{n,
      j} \bar{z} ) g_{n,j} ( z )|d|z|.
      \]
      Since $|g_{n,j}(z)|\le (\sum_j |g_{n,j}(z)|^q)^{1/q} <C$ and $\int|a_n
      (z'_{n,j} \bar{z} )|d|z|\le n^{-1}$, then 
      \[
      \|p\|_1\lesssim n^{-1}\sum_{j = 0}^{m_n} |p ( z'_{n, j} )|.
      \]
      \end{proof}

  To prove the inequality \eqref{necessaria} we will use the scheme proposed by
  Ramanathan and Steger in the context of the windowed Fourier transform (see
  \cite{RamSte95}). This works well when $p = 2$, for other $p \in [ 1, 
  \infty ]$ we use Lemma~\ref{p-estabilitat}. Now if we can prove the result
  for $p = 2$ we obtain the inequality
  \[ D^- ( \mathcal{Z} ) = \liminf_{R \to \infty} \left( \liminf_{n \to
     \infty} \frac{\min_{x \in [0,2\pi]} \# \mathcal{Z} ( n ) \cap ( x, x +
     R / n )}{R} \right) \ge \frac 1{2 \pi} - \delta . \]
  and this proves \eqref{necessaria} by taking $\delta$ arbitrarily small.
  
  Observe that the polynomial $p_n ( z ) = (z^n - 1) / ( 1 - z )$ has the
  property that
  \[ \int_{|z - 1| > R / n, |z| = 1} |p_n ( z ) |^2 \lesssim \frac{1}{R}
     \int_{|z| = 1} |p_n ( z ) |^2 . \]
  That means that for any separated family $\ZZ$ we have
  \begin{equation}
    \label{decay} \sum_{|z_{n, i} - 1| > R / n} |p_n ( z_{n, i} ) |^2 \lesssim
    \frac{1}{R} \int_{|z| = 1} |p_n ( z ) |^2 .
  \end{equation}
  Assume that $\ZZ$ is a $L^2$ M-Z family. Consider $\mathcal{P}_n$ the
  polynomials of degree $n$ as a Hilbert space with reproducing kernel. The
  corresponding reproducing kernel is $k ( z, w ) = ( 1 - ( z \bar{w} )^{n +
  1} ) / ( 1 - ( z \bar{w} ) )$, that is
  \[ p ( w ) = \langle p, k ( \cdot, w ) \rangle = \frac{1}{2 \pi} \int_{|z| =
     1} p ( z ) \overline{k ( z, w )} \, |d  z |,
     \quad \forall p \in \mathcal{P}_n . \]
  Since $\ZZ$ is a M-Z family that means that the normalized
  reproducing kernels $\{ \frac{1}{\sqrt{n}} k ( z, z_{n, i} ) \}_i$ form a
  frame in $\mathcal{P}_n$, i.e.
  \[ \| p \|^2 \simeq \frac{1}{n} \sum_{i = 1}^{m_n} | \langle p, k ( \cdot,
     z_{n, i} ) \rangle |^2, \quad \forall p \in \mathcal{P}_n \]
  with constants independent of $n$. This implies (see \cite{Daubechies92}
  for the basic facts on frames), that there are polynomials $\{ d_i ( z )
  \}_{i = 1}^{m_n}$ (the dual frame) such that for all polynomials $p$ in
  $\mathcal{P}_n$,
  \[ 
  \begin{split}
  p ( z ) = \frac{1}{\sqrt{n}} \sum_{i = 1}^{m_n} \langle p, k (
     z, z_i ) \rangle d_i ( z ), \\
     p ( z ) = \frac{1}{\sqrt{n}} \sum_{i = 1}^{m_n} \langle p, d_i ( z )
     \rangle k ( z, z_i ),
  \end{split} 
  \]
  and
  \[ \| p \|^2 \simeq \frac{1}{n} \sum_{i = 1}^{m_n} | \langle p, k ( \cdot,
     z_{n, i} ) \rangle |^2 \simeq \sum_{i = 1}^{m_n} | \langle p, d_i \rangle
     |^2, \quad \forall p \in \mathcal{P}_n \]
  Given $x \in \T$ and $t, r > 0$ ($t$ much bigger that $r$) we denote by $I (
  \tau )$ the arc-interval in $\T$ with center $x$ and radius $\tau / n$.
  consider the following two subspaces of $\mathcal{P}_n$:
  \begin{eqnarray*}
    W_S & = & \langle d_i ( z ) : z_i \in \ZZ ( n ) \cap I ( t + r ) \rangle\\
    W_I & = & \langle \frac{1}{\sqrt{n}} k ( z, w_j ) : w_j \in I ( t ), w_j^n
    = 1 \rangle .
  \end{eqnarray*}
  Let $P_S$ and $P_I$ denote the orthogonal projections of $\mathcal{P}_n$ on
  $W_S$ and $W_I$ respectively. We estimate the trace of the operator $T = P_I
  P_S$ in two different ways. To begin with
  \begin{equation}
    \label{tr1} \operatorname{tr} ( T ) \leq \operatorname{rank} W_S \leq
    \# ( \ZZ ( n ) \cap I ( t + r ) ) .
  \end{equation}
  On the other hand
  \[ \operatorname{tr} ( T ) = \sum_{w_i \in I ( t )} \langle T (
     \frac{1}{\sqrt{n}}
     k ( z, w_j ) ), P_I \kappa_j \rangle, \]
  where $\{ \kappa_j ( z ) \}$ is the dual basis of $\frac{1}{\sqrt{n}} k ( z,
  w_j )$ in $\mathcal{P}_n$. Using that $P_I$ and $P_S$ are projections one
  deduces that
  \begin{equation}
    \label{tr2} \operatorname{tr} ( T ) \geq \# \{ w_j \in I ( t ) \} \left( 1 -
    \sup_j | \langle P_S ( \frac{1}{\sqrt{n}} k ( z, w_j ) ) -
    \frac{1}{\sqrt{n}} k ( z, w_j ) ), \kappa_j \rangle | \right) .
  \end{equation}
  Since $\| \frac{1}{\sqrt{n}} k ( z, w_j ) \| \simeq 1$, also $\| \kappa_j \|
  \simeq 1$. We now show that $\| P_S ( \frac{1}{\sqrt{n}} k ( z, w_j ) ) -
  \frac{1}{\sqrt{n}} k ( z, w_j ) \| \leq \varepsilon$ for a suitable $r$.
  
  We have
  \[ \left| P_S ( \frac{1}{\sqrt{n}} k ( z, w_j ) ) - \frac{1}{\sqrt{n}} k (
     z, w_j ) \right|^2 \lesssim \frac{1}{n} \sum_{z_s \notin I ( t + r )} |
     \langle k ( z, w_j ), k ( z, z_s ) \rangle |^2 = \frac{1}{n} \sum_{z_s
     \notin I ( t + r )} \! \! |k ( w_j, z_s ) |^2 .
  \]
  This last sum is smaller than $\varepsilon$ if $r$ is big enough because
  $|w_j - z_s | \ge r / n$ and we can apply \eqref{decay}. If we put together
  \eqref{tr1} and \eqref{tr2}, we find that for every $\varepsilon$
there
  is an $r$ such that
  \[ \# ( \ZZ ( n ) \cap I ( t + r ) ) \ge ( 1 - \varepsilon ) \# \{ w_j \in
     I ( t ) \} = ( 1 - \varepsilon ) t, \]
  and this implies \eqref{necessaria}.
\end{proof}

The inequality in \eqref{necessaria} can be improved when $p = \infty$  to get
a strict inequality, thus providing a description in terms of
densities of the M-Z inequalities in this case. For this, we need to adapt
part of the arguments of Beurling in \cite{Beurling89b}. We will prove

\begin{theorem}\label{six}
  Let $p = \infty$. Given a separated family $\mathcal{Z}$ it is a M-Z family if
and only if
  \[ D^- ( \mathcal{Z} ) = \liminf_{R \to \infty} \left( \liminf_{n \to
     \infty} \frac{\min_{x \in [0,2\pi]} \# \mathcal{Z} ( n ) \cap ( x, x +
     R / n )}{ R} \right) > \frac 1{2\pi}, \]
\end{theorem}

\begin{definition}
  The Hausdorff distance between two compact sets $K, F$ in a metric space is
  defined as the infimum of the $\varepsilon > 0$ such that
  \[ K \subset ( F + B ( 0, \varepsilon ) ) \text{ and } F \subset ( K + B ( 0,
     \varepsilon ) ) . \]
  We denote this distance by $d_H ( K, F )$.
  
  A sequence of uniformly separated real sequences $\Lambda_n$ is said to
  converge weakly to $\Lambda$ if for any closed interval $I$, $d_H \left( ( I
  \cap \Lambda_n ) \cup \partial I, ( I \cap \Lambda ) \cup \partial I \right)
  \to 0$.
\end{definition}

\begin{definition}
  Recall that for any triangular family $\ZZ$ we can associate a sequence of
  real sequences $\Lambda ( n )$ as in \eqref{aplanats}. We take now an
  arbitrary family of real numbers $\tau_n$ and consider the corresponding
  translated sequences: $\Sigma ( n ) = \Lambda ( n ) - \tau_n$ (this
  corresponds to making rotations of $\ZZ ( n )$). We say that $\Lambda$
  belongs to a $W ( \ZZ )$ if there is a sequence of translates $\tau_n$ such
  that the corresponding $\Sigma ( n )$ converges weakly to $\Lambda$.
\end{definition}

\begin{definition}
  We denote by $\mathcal{F}$ the closed subspace of entire functions in the
  Bernstein class spanned by finite linear combinations of
  exponentials of the form $e^{irz}$ and $r \in \Q \cap [ - \pi, \pi ]$. The
  space $\mathcal{F}$ consists of almost periodic functions when restricted to
  the real line.
\end{definition}
With the same arguments as in \cite{Beurling89b} we can prove the following
theorem and corollary

\begin{theorem}
  The triangular family $\ZZ$ is a $L^{\infty}$ Marcinkiewicz-Zygmund family
  if and only if all $\Lambda \in W ( \ZZ )$ are uniqueness sets for
  $\mathcal{F}$.
\end{theorem}

\begin{corollary}
  If $\ZZ$ is a M-Z triangular family then there is an $\varepsilon > 0$ such
  that the triangular family $\ZZ'$ defined as $\ZZ' ( n ) = \ZZ ( [ n ( 1 -
  \varepsilon ) ] )$ is also a M-Z triangular family.
\end{corollary}
Now we apply the necessary condition \eqref{necessaria} of Theorem~\ref{tfive}
and we obtain that $D^- ( \ZZ)>D^- ( \ZZ' )\ge 2\pi $

\section{The model space} Actually it is possible to give a 
full characterization of M-Z sequences when $p=2$. It is not easily
computable. 
In this section we present this characterization. We need to introduce the
model spaces. Suppose that $I$
is an inner
function in the disk. We denote by
\[ K_I^2 ( \T ) = H^2 ( \T )  \ominus I  H^2 ( \T ) \]
If instead of the disk one considers the upper half plane, then $K_I^2 ( \R )$
is the standard $L^2$-Paley-Wiener space if $I = e^{iz}$. If we return back to
the disk and consider the case $I = z^n$ then $K_I^2$ is the space of
holomorphic polynomials of degree smaller or equal than $n$.

Thus, the setting of the model spaces is common for both the polynomials and
the Paley-Wiener space. Therefore any results that can be obtained from
general theorems in the model space setting will have the same flavor in
both the finite and the infinite-dimensional space.

Let us state the result that is more relevant in our
context. A Blaschke sequence $\Gamma \subset \D$ is a sampling sequence
for $K_I^2$ when
\[ \| f \|^2 \simeq \sum_{\Gamma} |f ( \gamma ) |^2 \omega_I ( \gamma ) , \]
for some appropriate weight $\omega_I$.
The following theorem was proved by Seip in \cite{Seip04}:
\begin{theorem}\label{seip} Denote by $B_\Gamma$ the Blaschke product with zeros
in $\Gamma$, 
  If $\Gamma$ satisfies $\sup_{\Gamma} |I ( \gamma ) |
  < 1$ and it is a Carleson sequence the following are equivalent:
  \begin{itemize}
    \item $\Gamma$ is a sampling sequence for $K_I^2$.
    
    \item There is an inner function $J$ such that 
    the Toeplitz operator in $H^2$ with symbol $JI \bar{B_\Gamma}$ is
    invertible.
    
  \end{itemize}
\end{theorem}

In our setting we start by a separated triangular family $\ZZ\subset \T$ and we
want a description of whether it is M-Z or not. We can replace this family by
the family $\WW$ defined as
\[
w_{n,j}= z_{n,j} (1-\epsilon/n)\qquad \forall j=0,\ldots,m_n,\ n\in\N. 
\] 
If $\epsilon>0$ is small enough the new triangular family is still separated
and by Lemma~\ref{perturbatiu} it will be a M-Z family whenever $\ZZ$ is a
M-Z
family. The advantage of $\WW$ is that we are uniformly under the hypotheses of
Theorem~\ref{seip}. That is if $I_n=z^n$ and $\Gamma_n$ is the sequence
$\WW(n)=\{w_{n,0},\ldots,w_{n,m_n}\}$, then $\sup_n\sup_{\Gamma_n}
|I_n(\gamma)|<1$
and moreover
$\Gamma_n$ is a Carleson measure (uniformly in $n$). Thus if define $B_n$ to be
the
Blashcke product with zeros in $\Gamma_n$, then a necessary and sufficient
condition so that $\ZZ$ is M-Z is that there exist inner functions $J_n$ such
that the Toeplitz operators $T_n$ in $H^2$ with symbols $J_n I_n \bar{B_n}$ are
invertible with uniform bounds. There are computable criteria for a Toeplitz
operator to be invertible (the Widom-Devinatz Theorem). The difficulty of
translating Theorem~\ref{seip} into a computable criteria are the inner
functions $J_n$. 
If we are given a sequence $\Gamma$ which we want to check wheter it is sampling
or not,  we do not have a natural candidate for function $J$ to
use the theorem. There are some
instances, for example in the Paley-Wiener 
space (see \cite{OrtSei02} or \cite{Seip04}) and for certain choices of
sequences $\Gamma$ where this is doable. In the finite
dimensional situation that we are dealing with, we do not get any new
computable criteria from this more complete theorem.


\def\cprime{$'$}
\providecommand{\bysame}{\leavevmode\hbox to3em{\hrulefill}\thinspace}
\providecommand{\MR}{\relax\ifhmode\unskip\space\fi MR }
\providecommand{\MRhref}[2]{%
  \href{http://www.ams.org/mathscinet-getitem?mr=#1}{#2}
}
\providecommand{\href}[2]{#2}

\end{document}